\documentclass[12pt]{amsart}
\usepackage{amscd,amssymb,graphicx}
\author{Alice Fialowski}
\address{E\"otv\"os Lor\'and University\\
Budapest, Hungary} \email{fialowsk@cs.elte.hu}
\author{Michael Penkava}
\address{University of Wisconsin\\
Eau Claire, WI 54702-4004} \email{penkavmr@uwec.edu}
\thanks{Research of the authors was partially supported by
grants from the University of Wisconsin-Eau Claire.}

\theoremstyle{definition}

\def \ph{\varphi}

\def \diag{\operatorname {diag}}

\def \refeq#1{equation (\ref{#1})}
\def \ra{\rightarrow}

\def \etc{\hbox{\it etc. }}

\def \tns{\otimes}

\def \mcom{,\cdots,}

\def \C{\mbox{$\mathbb C$}}

\def\d{d}

\def\linf{\mbox{$L_\infty$}}
\def\and{\mbox{ \rm and }}

\def\psa#1#2{\psi^{#1}_{#2}}

\def\P{\mathbb P}

\def\GL{{\mathbf{GL}}}

\def\sl{\mathfrak{sl}}

\begin{document}
\setlength{\multlinegap}{0pt}
\title{The moduli space of complex 5-dimensional Lie algebras}%

\date{\today}
\begin{abstract}
In this paper, we study the moduli space of all complex
5-dimensional Lie algebras, realizing it as a stratification by
orbifolds, which are connected only by jump deformations. The
orbifolds are given by the action of finite groups on very simple
complex manifolds. Our method of determining the stratification is
based on the construction of versal deformations of the Lie
algebras, which allow us to identify natural neighborhoods of the
elements in the moduli space.
\end{abstract}
\maketitle

\section{Introduction}
Moduli spaces of low dimensional Lie algebras possess very interesting properties. A moduli space of algebras of a fixed dimension
has a natural stratification by orbifolds. In the case of complex Lie algebras, these orbifolds consist of quasi-projective spaces,
given by removing a divisor from a certain $\C\P^n$, together with an action of a symmetric group. Some of the strata are singletons,
while other strata can be labelled by projective coordinates. This point of view is new, and uses deformation theory to give the
stratification in a unique form, where elements in a stratum deform in a smooth manner along the stratum, jump to elements in different
strata, and also deform smoothly along the families of the points to which they jump.

For algebras of dimension less than or equal to 5,
we obtained that each family (non singleton) of solvable Lie algebras contains one special element, which is nilpotent.  In fact, the
usual picture is that a family is given by the action of a symmetric group on $\P^n$, and the \emph{generic point} of $\P^n$ corresponds to
the nilpotent.  We remark that if $\P^n$ is given by projective coordinates $(p_0:\cdots:p_n)$, then the generic point is $(0:\cdots:0)$,
which is usually excluded from consideration by algebraic geometers, but cannot be excluded in our consideration, since there is an algebra
corresponding to the generic point. (The terminology \emph{generic point} is a bit unfortunate, since the generic point behaves in a very
nongeneric manner.) The key point here is that the nilpotent elements fit within a larger picture in a natural manner.  We also note that
the generic element in one family may be isomorphic to the generic element in another family.  This behavior represents the only overlap
between families.

A big difference between our analysis of the structure of the moduli space of algebras, and classifications given by previous authors
is that our goal is not to simply determine a list of the algebras, but to see how the moduli space is glued together. For that, we needed
to compute not only the cohomology of the algebras, but \emph{versal deformations} of the algebras, which give complete information about all
of the deformations. We will describe later how these versal deformations are computed.

A classification of two and three dimensional Lie algebras is given in \cite{jac}. A classification
of real four and five dimensional algebras was given in \cite{mub1,mub2}.  A list of the 5-dimensional Real Lie algebras
was given in \cite{PSWZ} and the same list also appears in the recent article {FGH}. The main emphasis of the classification methods
was to give a complete list of the algebras, without duplicates, although sometimes the lists were subdivided in terms
of other properties, such as invariants of the algebras.

The description we give of the moduli space gives a more natural decomposition of the moduli space, and divides it into fewer pieces.
We first give a summary of the moduli spaces of  three and four dimensional complex Lie algebras, then we give an analysis of how we
construct the moduli space of 5-dimensional algebras, and afterwards, we give details on the deformations of the algebras.

\section{Three Dimensional Lie algebras}
In \cite{fp3,Ott-Pen}, a decomposition of the moduli space of 3-dimensional complex Lie algebras was given, in terms of strata given by projective orbifolds.
However, in \cite{Ott-Pen}, where the projective description was given explicitly, we still were missing one piece of the puzzle. In the space $\C\P^1$ (hereafter referred to as $\P^1$), there
is a special element $(0:0)$, called the \emph{generic element}, which is usually excluded from the space, because with the removal of this element
$\P^1$ becomes a smooth manifold, but when this element is included, the topology becomes not even Hausdorff.  However, in the correspondence between the elements of $\P^1$ to the algebras labeled as $d_2(p:q)$, it turns out that setting both $p$ and $q$ equal to zero gives rise to an algebra, and the manner
in which this algebra deforms in relation to the family $d_2(p: q)$ corresponds in a natural manner to the way in which the generic element sits in the
space $\P^1$. This means that there are jump deformations from the algebra $d_2(0:0)$ to every element in the family, in parallel with the fact that every open nbd of the point $(0: 0)$ in $\P^1$ includes all of the points in $\P^1$.  Thus, the final form of the description of the moduli space, representing a completely canonical ordering in terms of deformation theory, is given in the table below.
\begin{table}[ht]
\begin{center}
\begin{tabular}{llrrrrr}
Type&Formula&$H^0$&$H^1$&$H^2$&$H^3$\\
\hline
$d_1$&$\psa{12}3+\psa{13}2+\psa{23}1$&0&0&0&0\\
$d_2(p:q)$&$\psa{13}1p+\psa{23}1+\psa{23}2q$&0&1&1&0\\
$d_2(1: 1)$&$\psa{13}1+\psa{23}1+\psa{23}2$&0&1&2&1\\
$d_2(1:0)$&$\psa{13}1+\psa{23}1$&1&2&1&0\\
$d_2(0:0)$&$\psa{23}1$&1&4&5&2\\
$d_3$&$\psa{13}1+\psa{23}2$&0&3&3&0\\
\hline\\
\end{tabular}
\end{center}
\caption{\protect Cohomology of Three Dimensional Complex Lie
Algebras}\label{Table 1}
\end{table}

Although the algebra given as $d_2(p:q)$ is labelled by $\P^1$, there is an identification of elements $d_2(p:q)$ and $d_2(q:p)$, that is to say,
these two elements are isomorphic.  This means that $d_2(p:q)$ is labelled by the projective orbifold $\P^1/\Sigma_2$, where the
symmetric group $\Sigma_2$ acts on $\P^1$ by permuting the coordinates.

The algebra listed as $d_1$ is $\sl_2(\C)$, and it occurs first in the list because it is rigid, and so does not deform to any other algebra. The
family $d_2(p:q)$ roughly corresponds to the algebras which are classically referred to as $\mathfrak{r}_3(\lambda)$, where $\lambda=q/p$, except for the special cases
$d_2(1:0)=\mathfrak{r}_2(\C)\oplus\C$, $d_2(1:1)=\mathfrak{r}_3(\C)$, and $d_2(0:0)=\mathfrak{n}_3$. Finally, $\mathfrak{r}_3(1)=d_3$.

Now, the most serious change in our notation is that the elements $\mathfrak{r}_3(1)$ and $\mathfrak{r}_3(\C)$ have been interchanged, that is, we have altered the family by removing one of the elements and replacing it with another.  This is a very important change, which has to do with the natural projective decomposition of the moduli space of $2\times 2$ matrices with the action of the group $\mathbf{SL}_2(\C)\times\C^*$ given by conjugation by the matrix and multiplication by the element in $\C^*$. The equivalence classes of matrices correspond to the Jordan decompositions of the matrices, up to a multiplication of
the diagonal elements by a constant.  The classical decomposition corresponds to the following decomposition
\begin{equation*}
\mathfrak{r}_3(\lambda)\leftrightarrow\left[\begin{array}{rr}1&0\\0&\lambda\end{array}\right],
\quad\mathfrak{r}_3(\C)\leftrightarrow\left[\begin{array}{rr}1&1\\0&1\end{array}\right],
\quad\mathfrak{n}_3(\C)\leftrightarrow\left[\begin{array}{rr}0&1\\0&0\end{array}\right].
\end{equation*}
Our decomposition corresponds to
\begin{equation*}
d_2(p:q)\leftrightarrow\left[\begin{array}{rr}p&1\\0&q\end{array}\right],\qquad
d_3\leftrightarrow\left[\begin{array}{rr}1&0\\0&1\end{array}\right].
\end{equation*}
The key difference is that we include the elements whose Jordan form has a single block in the family, and this is for a very important reason.
If we make a small change in the identity matrix, adding a small constant $t$ in upper right corner, to obtain the matrix
$\begin{smallmatrix}1&t\\0&1\end{smallmatrix}$, it is easy to see that this matrix is equivalent to the matrix
$\begin{smallmatrix}1&1\\0&1\end{smallmatrix}$. Moreover, this is not only true for the matrices, but for the algebras they correspond to, so there
is a jump deformation from $d_3$ to $d_2(1: 1)$, but not the other way around. Accordingly, the correct way to decompose the moduli space in terms of
the ideas of deformation theory is in terms of the decomposition we have given, rather than the classical decomposition.

We have not yet given a description of how a $2\times 2$ matrix is related to a 3-dimensional algebra. More generally, if a Lie algebra can be
decomposed as an extension of a 1-dimensional Lie algebra by a trivial $n$-dimensional algebra, then the algebra is completely determined by the
module structure on the $n$-dimensional space as a representation of the trivial algebra. In terms of the codifferential notation, we have
that the algebra is of the form $d=\psi^{j,n+1}_ia^i_j$, where the matrix $A=(a^i_j)$ is any $n\times n$ matrix.  Thus the algebras are determined by an $n\times n$ matrix, and any such matrix determines a Lie algebra. Moreover, the equivalence classes of algebras of this type are precisely given by the equivalence classes of Jordan decompositions, up to a scaling factor of the eigenvalues. Understanding
the deformations of the algebras within this class of algebras is equivalent to understanding how the matrices can deform into each other, which gives
a decomposition into strata in a precise manner. This Jordan decomposition gives the family $d_2(p:q)$, and $d_3$, while the zero matrix corresponds
to the trivial Lie algebra structure, which we omit from the description of the moduli space.

Note that the algebra $d_2(0:0)=\mathfrak{n}_3$ fits into this family, and it has jump deformations to all the members of the family, as well as to the
algebra $d_1=\mathfrak{sl}_2(\C)$. But it does not jump to the element $d_3=\mathfrak{r}_3(1)$. This is another vindication of the classification scheme we
have introduced for the 3-dimensional complex Lie algebras.

\section{Four Dimensional Complex Lie Algebras}
The classification of 4-dimensional complex Lie algebras was first given in \cite{mub2}. Actually, the author gave a classification of real Lie algebras in this paper.  Several others have given classifications of complex Lie algebras, for example, \cite{AG,BS}. In \cite{fp8}, we gave a decomposition of this moduli space into
projective orbifolds, corresponding to the deformation theory of the algebras.  Our notation in this paper was somewhat eclectic, and we had not yet
realized the significance of the \emph{generic element} in projective space. It turns out that the moduli space decomposes into exactly 7 families (some singletons).  The table below gives the cohomology of every element in the moduli space, including special cases within families where the
cohomology does not follow the generic pattern.

\begin{table}[ht]
\begin{center}
\begin{tabular}{lrrrrr}
Type&$H^0$&$H^1$&$H^2$&$H^3$&$H^4$\\
\hline
$d_1=\psa{23}4+\psa{24}3+\psa{34}2$&1&1&0&1&1\\
$d_2=\psa{13}1+\psa{24}2$&0&0&0&0&0\\\hline
$d_3(p:q)=\psa{23}1+\psa{14}1(p+q)+\psa{24}2p+\psa{34}2+\psa{34}3q$&0&1&1&0&0\\
$d_3(1:0)$&0&1&2&1&0\\
$d_3(1:-1)$&1&2&2&2&1\\
$d_3(0:0)$&1&4&6&5&2\\
$d_4=\psa{23}1+\psa{14}12+\psa{24}2+\psa{34}3$&0&3&3&0&0\\\hline
$d_5(p:q:r)=\psa{14}1p+\psa{24}1+\psa{24}2q+\psa{34}2+\psa{34}3r$&0&2&2&0&0\\
$d_5(p:q:-p-q)$&0&2&2&1&1\\
$d_5(p:q:p+q)$&0&2&3&1&0\\
$d_5(p:q:0)$&1&3&3&1&0\\
$d_5(1:-1:0)$&1&3&5&5&2\\
$d_5(0:0:0)$&1&4&6&5&2\\
\hline
$d_6(p:q)=\psa{14}1p+\psa{24}2p+\psa{34}2+\psa{34}3q$&0&4&4&0&0\\
$d_6(1:-2)$&0&4&4&1&1\\
$d_6(1:2)$&0&4&5&1&0\\
$d_6(0:1)$&2&6&6&2&0\\
$d_6(1:0)$&1&5&7&3&0\\
$d_6(0:0)$&2&8&13&10&3\\
$d_7=\psa{14}1+\psa{24}2+\psa{34}3$&0&8&8&0&0\\
\hline\\
\end{tabular}
\end{center}
\caption{\protect Cohomology of Four Dimensional Complex Lie
Algebras}\label{Table 2}
\end{table}

The algebra $d_3(p:q)$ is labelled by $\P^1/\Sigma_2$, where $\Sigma_2$ acts in the usual manner by permuting the projective coordinates $(p:q)$.
The algebra $d_5(p:q:r)$ is labelled by $\P^2/\Sigma_3$, where the symmetric group $\Sigma_3$ acts on $\P^2$ by permuting the coordinates.
The algebra $d_6(p:q)$ is labelled by $\P^1$, without any action of the symmetric group.  We say that these algebras determine strata given by
projective orbifolds, which are given by an action of a group on a projective space $\P^n$.

The algebras  $d_5(p:q:r)$, $d_6(p:q)$ and $d_7$ arise as extensions of the 1-dimensional Lie algebra by the trivial 3-dimensional Lie algebra, so are determined by a $3\times 3$ matrix representing the module structure. We have
\begin{equation*}
d_5(p:q:r)\leftrightarrow\left[\begin{array}{rrr}p&1&0\\0&q&1\\0&0&r\end{array}\right],
d_6(p:q)\leftrightarrow\left[\begin{array}{rrr}p&0&0\\0&p&1\\0&0&q\end{array}\right],
d_7\leftrightarrow\left[\begin{array}{rrr}1&0&0\\0&1&0\\0&0&1\end{array}\right].
\end{equation*}
More precisely, we mean that the algebras correspond to these matrices, which give the Jordan decomposition of the types of $3\times 3$ matrices up to scaling the diagonal elements.

The algebras $d_3(p:q)$ and $d_4$ correspond to extensions of the trivial 1-dimensional algebra by the only nontrivial nilpotent 3-dimensional algebra
$\mathfrak{n}_3$. The matrices representing the module structure on these algebras are
\begin{equation*}
\d_3(p:q)\leftrightarrow\left[\begin{array}{rrr}p+q&0&0\\0&p&1\\0&0&q\end{array}\right],\qquad
d_4\leftrightarrow\left[\begin{array}{rrr}2&0&0\\0&1&0\\0&0&1\end{array}\right].
\end{equation*}
One sees a similar Jordan decomposition of the lower right $2\times 2$ submatrix.

Now the generic elements $d_3(0:0)$, $d_5(0:0:0)$ and $d_6(0:0)$ are nilpotent.  However, there are only 2 nontrivial nilpotent complex 4-dimensional Lie algebras, which seems contradictory, until one realizes that $d_3(0:0)$ and $d_5(0:0:0)$ are isomorphic.

Every finite dimensional, solvable, complex
Lie algebra $L$ has a nilpotent ideal whose dimension is at least half of the dimension of the algebra. In fact, there is a unique maximal nilpotent ideal $\mathfrak n$, called the nilradical, which includes the derived algebra $[L,L]$. This makes it much easier to construct all the algebras.  The algebra $d_1$ is the direct sum of $\mathfrak{sl}_2(\C)$ and the trivial algebra $\C$, so is not solvable. The algebra
$d_2$ is $\mathfrak{r}_2(\C)\oplus\mathfrak{r}_2(\C)$, a direct sum of two copies of the nontrivial (solvable) 2-dimensional complex Lie algebra.
However, it also can be realized as an extension of the trivial 2-dimensional algebra by the trivial 2-dimensional algebra, so it has
a nilpotent ideal of dimension 2, and in fact, this ideal is the nilradical.

The algebras $d_3(p:q)$ and $d_4$ are extensions of the trivial 1-dimensional Lie algebra by the nontrivial 3-dimensional nilpotent algebra.
The algebras $d_5(p:q:r)$, $d_6(p:q)$ and $d_7$ are extensions of the trivial 1-dimensional Lie algebra by the trivial 3-dimensional Lie algebra.

It is very interesting that every nontrivial 4-dimensional nilpotent algebra arises as the \emph{generic element} in a family. This would lead
one to suspect that
this pattern might be a general pattern, although the existence of families of nilpotent Lie algebras, the filiform families, in higher dimensions, means that the pattern is more complicated in general.
It is also clear that there can be families of Lie algebras whose generic element is not nilpotent, if we consider
a direct sum of a simple Lie algebra and a family. As we shall see, in dimension 5, there is a nilpotent Lie algebra which is a singleton, not a member
of a family.

So far, for dimensions $3$ and $4$, we have found a natural decomposition of the moduli space into strata, each of which is a projective orbifold,
and moreover, the strata are ``glued together" by jump deformations. The fact that the moduli space can be decomposed into strata has been known
for a long time, but the description of the strata as projective orbifolds was a new discovery of the authors and some collaborators.

\section{Five Dimensional Lie Algebras}
If we have an exact sequence
$$0\ra M\ra L\ra W\ra 0$$
of algebras, we say that $L$ is an extension of $W$ by $M$. We can view $M$ as an ideal in $L$ and $W$ as the quotient algebra $L/M$.  Denote the algebra structures on
$M$, $L$, and $W$ by $\mu$, $d$, and $\delta$ resp. Then the algebra structure $d=\delta+\mu+\lambda+\psi$, where
$\lambda:M\tns W\ra M$ is the ``module structure" and $\psi:W\tns W\ra M$ is the ``cocycle". In terms of the brackets of \emph{cochains},
the $d$ determines an algebra structure precisely when the following three conditions are satisfied:
\begin{align*}
&[\mu,\lambda]=0, \text{the compatibility condition}\\
&\frac12[\delta+\lambda,\delta+\lambda]+[\mu,\psi]=0, \text{the Maurer-Cartan condition}\\
&[\delta+\lambda,\psi]=0,\text{the cocycle condition}.
\end{align*}
When $\mu=0$, the Maurer-Cartan condition (MC condition) reduces to $[\delta+\lambda,\delta+\lambda]=0$, which is the classical MC condition, and it is only when this classical condition holds that $\lambda$ determines a module structure on $M$ over the algebra structure $\delta$ on $W$,  and in
that case, $\psi$ is actually a cocycle with respect to a coboundary operator $D_{\delta+\lambda}$ given by $D_{\delta+\lambda}(\ph)
=[\delta+\lambda,\ph]$. The MC condition guarantees that $D^2_{\delta+\lambda}=0$. We retain the names module structure for $\lambda$ and cocycle
for $\delta$ in general, although they are not strictly speaking accurate.

We first classify the nonsolvable 5-dimensional Lie algebras, then turn our attention to the solvable ones.

\subsection{Nonsolvable 5-dimensional Lie algebras}
If a Lie algebra $L$ is not solvable, then it has a decomposition
\begin{equation*} 0\ra M\ra L\ra S\ra 0,\end{equation*}
where $S$ is a semisimple algebra, and $M$ is the maximal solvable ideal in $L$.
When the Lie algebra is defined over $\C$, one can say a bit more, because we have a \emph{Levi decomposition} $L=M\rtimes S$, where $M$ is the maximal solvable ideal and $S$ is a semisimple subalgebra of $L$. This means that in the extension of the algebra structure $\delta$ on $S$
by the algebra structure $\mu$ on $M$, we have $\psi=0$, so we only have to consider the module structure $\lambda$, which actually will
be a module structure, since the term $[\mu,\psi]$ vanishes, so the classical MC condition holds.

Since the only semisimple complex Lie algebra of dimension less than or equal to $5$ is $S=\sl(2,\C)$, we know that a non-solvable Lie algebra
is an extension of $S$ by a 2-dimensional solvable algebra, which can either be the unique nontrivial 2-dimensional Lie algebra, or  the the trivial one.  Suppose $M=\langle e_1,e_2\rangle$ and $S=\langle e_3,e_4,e_5\rangle$, by which we mean that $M$ has the ordered basis $e_1,e_2$ and similarly for $W$. Then we can represent $\sl(2,\C)$ by the codifferential
$\delta=2\psa{34}4-2\psa{35}5+\psa{34}5$.

The nontrivial solvable Lie algebra of dimension 2 can be represented by the codifferential $\mu=\psa{12}1$. The only extension of $\delta$ by $\mu$
is the direct sum which gives the first algebra on our list: $$d_1=\psa{12}1+2\psa{34}4-2\psa{35}5+\psa{34}5.$$

For extensions by the trivial structure, these correspond to the $\sl(2,\C)$ module structures on $\C^2$. There is a standard irreducible representation of $\sl(2,\C)$ on $\C^2$, which gives the algebra
$$d_2=2\psa{34}4-2\psa{35}3+\psa{45}3-\psa{13}1+\psa{2,3}2+\psa{24}1+\psa{15}2.$$
The other extension is just the direct sum of $\delta$ and the trivial algebra structure on $M$, so is given by
$$d_3=2\psa{34}4-2\psa{35}5\psa{34}5.$$

This completes the classification of the nonsolvable 5-dimensional Lie algebras.

\subsection{Solvable 5-dimensional Lie algebras}
The nilradical of a solvable complex Lie algebra has dimension at least 3, by a classical Lie theory result, and the nilradical $\mathfrak n$ contains the derived algebra $L'=[L,L]$, which means that $W=L/\mathfrak n$ is a trivial algebra.
Thus every 5 dimensional solvable Lie algebra can be given as an
extension of a trivial algebra of dimension 1 or 2 by a nilpotent algebra of complementary dimension 4 or 3. We begin by classifying the Lie algebras
whose nilradical is at least 4-dimensional.  There are three possibilities which we have to consider, $\mu=\psa{34}2$, $\mu=\psa{23}1+\psa{34}2$ and
the trivial algebra.  Then we will classify those algebras which have a nilradical of dimension 3, which can either be the algebra $d=\psa{23}1$ or
the trivial algebra.
\subsection{Extensions of a 1-dimensional algebra by a 4-dimensional trivial algebra}
There are 5 classes of $4\times 4$ matrices, essentially determined by the Jordan decomposition of the matrix (up to a constant multiple), and
these yield 5 algebras, 4 of which are families.  The matrices are as follows.
\begin{align*}
\left[\begin{matrix}p&1&0&0\\0&q&1&0\\0&0&r&1\\0&0&0&s\end{matrix}\right],
\left[\begin{matrix}p&0&0&0\\0&p&1&0\\0&0&q&1\\0&0&0&r\end{matrix}\right],
\left[\begin{matrix}p&1&0&0\\0&q&0&0\\0&0&p&1\\0&0&0&q\end{matrix}\right],
\left[\begin{matrix}p&0&0&0\\0&p&0&0\\0&0&p&1\\0&0&0&q\end{matrix}\right],
\left[\begin{matrix}1&0&0&0\\0&1&0&0\\0&0&1&0\\0&0&0&1\end{matrix}\right].
\end{align*}
The first matrix gives rise to the family $d_{20}(p:q:r:s)$ which is parameterized by $\P^3/\Sigma_4$, where the symmetric group $\Sigma_4$ acts
on $\P^4$ by permuting the projective coordinates.  The element
$d_{20}(0:0:0:0)=\psa{25}1+\psa{34}2+\psa{45}3$, the generic element in the family, is nilpotent.

The second matrix gives rise to the family $d_{21}(p:q:r)$ which is parameterized by $\P^2/\Sigma_2$, which acts by permuting the last two coordinates.  The generic element in the family, $d_{21}(0:0:0)=\psa{35}2+\psa{45}3$, is again a nilpotent algebra.

The third matrix gives rise to the family $d_{22}(p:q)$ which is parameterized by $\P^1/\Sigma_2$, which acts by permuting the coordinates.
The generic element $d_{22}(0:0)=\psa{25}1+\psa{45}3$ is also nilpotent.

The fourth matrix gives rise to the family $d_{23}(p:q)$, which is parameterized by $\P^1$ with no symmetries. The generic element $d_{23}(0:0)=\psa{45}3$ is nilpotent.

Finally, the fifth matrix gives rise to the algebra $d_{24}$.

The classification of the part of the moduli space given by extensions of the trivial 1-dimensional algebra by a trivial nilpotent algebra is very
easy, and it is also easy to order the elements by simply looking at the cohomology of a generic element in one of the subfamilies. The ideas
we presented here extend in a straightforward manner to describe the extensions of a 1-dimensional (necessarily trivial) algebra by an $n$-dimensional trivial algebra.

\subsection{Extensions of a 1-dimensional algebra by the 4-dimensional nilpotent algebra $\mu=\psa{34}2$}
For extensions of this type, we have to take into account that if $\lambda$ is the module structure on $M$, then the compatibility condition
$[\mu,\lambda]=0$ imposes a nontrivial constraint.  The matrix $A=(a^i_j)$ of $\lambda=\psa{j5}ia^i_j$ must be of the form
\begin{align*}
A=\left[\begin{matrix}
p&0&a_{13}&a_{14}\\
a_{21}&q+r&0&0\\
0&0&q&a_{34}\\
0&0&a_{43}&r
\end{matrix}\right]
\end{align*}

The group $G_{\delta,\mu}$ of linear transformations preserving $\mu$ and $\delta$  is the direct product of the group $G_\mu$ of linear transformations of $\C^4$
which preserve $\mu$ and the group $\C^*$, (which is the group of linear transformations of $\C$ preserving $\delta=0$)and the action
is given by conjugation of the matrix $A$ of $\lambda$ by the matrix $G$ of an element in $G_\mu$ up to multiplication by a number.
An element of this groups acts on the space of algebras of this form, and two algebras related by this action are isomorphic.  The
action allows us to determine an element in one of these equivalence classes which has a nice form, without losing any isomorphism types.

A matrix $G$ represents an element in $G_\mu$ precisely when it is of the form
$$G=\left[ \begin {array}{ccccc} g_{{1,1}}&0&g_{{1,3}}&g_{{1,4}}&0
\\\noalign{\medskip}g_{{2,1}}&g_{{4,4}}g_{{3,3}}-g_{{3,4}}g_{{4,3}}&g_
{{2,3}}&g_{{2,4}}&0\\\noalign{\medskip}0&0&g_{{3,3}}&g_{{3,4}}&0
\\\noalign{\medskip}0&0&g_{{4,3}}&g_{{4,4}}&0\\\noalign{\medskip}0&0&0
&0&g_{{5,5}}\end {array} \right],$$ and its determinant
$g_{{1,1}} \left( g_{{4,4}}g_{{3,3}}-g_{{3,4}}g_{{4,3}} \right) ^{2}g_{
{5,5}}$ is nonzero. Let us assume that $g_{21}$, $g_{14}$, $g_{24}$, $g_{23}$ and $g_{13}$ vanish, which will not affect the determinant.
Then the elements $a_{13}$ $a_{14}$ are replaced by $a_{13}g_{33}+a_{14}g_{43}$ and $a_{13}g_{34}+a_{14}g_{44}$, so this represents the
action of $\GL(2)$ on $\C^2$, which means that we can assume that $a_{14}=0$ and either $a_{13}=1$ or $a_{13}=0$.
Next, it turns out that if $a_{34}\ne0$ we can assume that $a_{43}=0$, so we get three subcases given by $a_{34}=1$ and $a_{43}=0$,
$a_{34}=0$ and $a_{43}=1$, and $a_{34}=a_{43}=0$.  Finally, for each of these cases, we can assume that either $a_{21}=1$ or $a_{21}=0$.
Thus we split the analysis into a bunch of discrete cases. Of course, these cases overlap quite a bit, but after some analysis, we obtain
the following families, including some singletons:
\begin{align*}
d_{12}(p:q:r)&=\psa{34}2+\psa{15}1p+\psa{15}2+\psa{25}2(q+r)+\psa{35}1+\psa{35}3q+\psa{45}3+\psa{45}4r\\
d_{13}(p:q)&=\psa{34}2+\psa{15}1(p+q)+\psa{25}2(p+q)+\psa{35}1+\psa{35}3p+\psa{45}3+\psa{45}4q\\
d_{14}(p:q)&=\psa{34}2+\psa{15}1p+\psa{15}2(p+q)+\psa{35}1+\psa{35}3q+\psa{35}4+\psa{45}4p\\
d_{15}(p:q)&=\psa{34}2+\psa{15}1p+\psa{15}2+\psa{25}2(2q)+\psa{35}1+\psa{35}3q+\psa{45}4q\\
d_{16}&=\psa{34}2+2\psa{15}1+2\psa{25}2+\psa{35}1+\psa{35}3+\psa{45}4\\
d_{17}&=\psa{34}2+\psa{15}1+\psa{15}2+2\psa{25}2+\psa{35}3+\psa{45}4\\
d_{18}&=\psa{34}2+\psa{15}1+\psa{25}2+\psa{35}1+\psa{35}4+\psa{45}4\\
d_{19}&=\psa{34}2+\psa{15}2
\end{align*}

The family $d_{12}(p:q:r)$ is parameterized by $\P^2/\Sigma_2$, where $\Sigma_2$ acts by permuting the last two coordinates.
The family $d_{13}(p:q)$ is parameterized by $\P^1/\Sigma_2$, where $\Sigma_2$ acts in the standard manner.
The families $d_{14}(p:q)$ and $d_{15}(p:q)$ are parameterized by $\P^1$, with no action of a symmetric group.

\subsection{Extensions of a 1-dimensional algebra by the 4-dimensional nilpotent algebra $\mu=\psa{23}1+\psa{34}2$}
The compatibility condition $[\mu,\lambda]=0$ places strong limitations on the form of the matrix
for $\lambda$. Moreover, if $\lambda$ is replaced by $\lambda+[\mu,\beta]$ where $\beta:W\ra M$ is a 1-cochain, then we arrive at
an isomorphic algebra, so we can use this to express $\lambda$ in a nice form.  Its $4\times 4$ matrix $A$ can be given by
$$A=\left[\begin {array}{cccc} 2\,a_{{3,3,1}}+a_{{4,4,1}}&0&0&a_{{1,4,1}}\\\noalign{\medskip}0&a_{{3,3,1}}+a_{{4,4,1}}&0&0\\\noalign{\medskip}0
&0&a_{{3,3,1}}&0\\\noalign{\medskip}0&0&a_{{4,3,1}}&a_{{4,4,1}}
\end {array} \right].$$ Moreover, considering the action of the group $G_{\mu,\delta}$ of symmetries of $\mu$, we can also assume that
$a_{14}$ and $a_{43}$ are either 1 or 0, which gives 4 cases. We obtain the three algebras
\begin{align*}
d_9(p:q)=&\psa{23}1+\psa{34}2+\psa{15}1(2p+q)+\psa{25}2(p+q)\\&+\psa{35}3p+\psa{35}4+\psa{45}1+\psa{45}4q\\
d_{10}=&\psa{23}1+\psa{34}2+3\psa{15}1+2\psa{25}2+\psa{35}3+\psa{45}1+\psa{45}4\\
d_{11}=&\psa{23}1+\psa{34}2+\psa{15}1+\psa{25}2+\psa{35}4+\psa{45}4\
\end{align*}

This completes the analysis of the algebras which have a nilradical of dimension 4.
We now turn our attention to the algebras with a nilradical of dimension 3. Their classification is more subtle than the ones we have
studied so far.

\subsection{Lie algebras with trivial 3-dimensional nilradical}
A Lie algebra with trivial nilradical $M=\langle e_1,e_2,e_3\rangle$ is determined by the module structure
$\lambda=\psa{j4}iA^i_j+\psa{j5}iB^i_j$, and the cocycle $\psi=\psa{45}ic^i$. The algebra structure on $W$ is $\delta=0$, the algebra
structure on $M$ is $\mu=0$  and the MC condition
reduces to $[\lambda,\lambda]=0$.  But this happens precisely when the matrices $A$ and $B$ in the definition of $\lambda$ commute.
Moreover, there is an action of $G_{M,W}=\GL_3\times\GL_2$ on the $\lambda s$, which is as follows: If $G=\diag{G_M,G_W}$ is a block diagonal
decomposition of an element of $G_{M,W}$, and $\lambda'=\psa{j4}i\hat A^i_j+\psa{j5}i\hat B^i_j$, then $\hat A$ and $\hat B$ are given by
linear combinations of the conjugation of $A$ and $B$ by $G_M$, where the coefficients of the linear combination are determined by $G_W$.

We first note that if one of the matrices $A$ or $B$ is strictly upper triangular, then the maximal nilpotent ideal will have dimension larger than
three, so we arrange for at least generically this not to be true.  A strictly upper triangular matrix can be conjugated to one of the following forms
\begin{equation}\label{threematrices}\left[\begin{matrix}p&1&0\\0&q&1\\0&0&r \end{matrix}\right],\qquad\left[\begin{matrix}p&0&0\\0&p&1\\0&0&q\end{matrix}\right],\qquad
\left[\begin{matrix}1&0&0\\0&1&0\\0&0&1\end{matrix}\right].
\end{equation}

Let us suppose that $A$ has been conjugated to the first of these forms.  Then the conjugated form of $B$ must be upper triangular.  Note that we can assume that neither $A$ nor $B$ has a zero diagonal, because in that case, the maximal nilpotent ideal of $L$ would have higher dimension. Thus
by adding a multiple of the conjugated $B$ to the conjugated $A$, and then restoring it to Jordan form, we can assume $A$ is of the form
$A=\left[\begin{matrix}0&1&0\\0&p&1\\0&0&q\end{matrix}\right]$, and $B$ is upper triangular. We can also assume $p\ne0$, because one of $p$ or $q$ can be assumed nonzero, so we can conjugate to assume that $p$ does not vanish. By adding a multiple of the conjugated matrix $A$ to the conjugated
matrix $B$, we can thus assume that $B^2_2=0$. After these assumptions, we can then express $B$ uniquely in the form
$B=\left[\begin{matrix}(-r+uq)p&r-uq&u\\0&0&r\\0&0&r(q-p)\end{matrix}\right]$.  Since scaling $e_5$ by a constant multiplies the matrix $B$ by an
overall constant, we can reduce to two cases, $u=1$ or $u=0$.

\subsubsection{The case $u=1$}
Let us first consider the case $u=1$. Then $B$ has the form $B=\left[\begin{matrix}(-r+q)p&r-q&1\\0&0&r\\0&0&r(q-p)\end{matrix}\right]$. The first
question we ask about this form is whether it is projective, that is, whether the algebra determined by $(p,q,r)$ is isomorphic to the algebra determined by $(tp,tq,tr)$, when $t\ne 0$, and it can be shown that this is true.

Now, it is true that one also has to consider the possible values of $\psi$, but it can be shown that whenever $p\ne$ and not both $q$ and $r$ vanish, that the algebra determined by adding $\psi$ is isomorphic to the algebra obtained by ignoring it.

Thus we obtain a projective family $d_5(p:q:r)$.  However, there are problems with this ``family'', because, for certain subfamilies, there are more isomorphisms than usual, so that in order to obtain a good family, we have to exclude 7 $\P^1$s, given by $p=0$, $q=0$, $r=0$, $p=q$, $p=r$,
$q=r$, and $pr=q^2$.
Then, on the resulting space $S$, which is $\P^2$ minus a certain divisor, there is an action of $\Sigma_3$, but not the usual one we have encountered, so that $d_5$ is given by $S/\Sigma_3$.  The group of symmetries is generated by the transformations
\begin{align*}
\tau(p: q: r)&={\frac {{q}^{2} \left( -q+p \right) }{rp-{q}^{2}}},-{\frac {{q}^{3} \left( r-q \right) }{r \left( rp-{q}^{2} \right) }},{\frac {p \left(
r-q \right) ^{2}{q}^{2}}{r \left( -q+p \right)  \left( rp-{q}^{2}
 \right) }}\\
\sigma(p:q:r)&=\left({\frac {r \left( p-q \right) ^{2}}{ \left( -q+r \right)  \left( r-p \right) }}:{\frac {p \left( p-q \right) }{r-p}}:{\frac { \left( r-q
 \right) p}{r-p}}\right).
\end{align*}
It can be checked that $\tau$ has order 2, $\sigma$ has order 3, and that $\sigma\tau=\tau\sigma^2$, so that they generate a group isomorphic  to
$\Sigma_3$.

Note that there are some necessary conditions on $(p:q:r)$ in order for these maps to be well defined. It is immediate from the formulas that
we cannot have $r=0$, $p=q$, $p=r$, $q=r$ or $pr=q^2$. However, the fact that $p\ne0$ and $q\ne0$ can be determined by examining powers of these
maps.


Now, let us analyze the seven $\P^1$s which we have excluded.  The first $\P^1$ is given by elements of the form $d_5(0:p:q)$.  It turns out that
$d_5(0:p:q)\sim d_5(p:0:0)\sim d_{15}(p:0)$.  The entire $\P^1$ collapses to 1 point (excluding the generic element), and this algebra is isomorphic
to a special case in another family. Thus this $\P^1$ doesn't generate any new algebras.

The second $\P^1$ consists of the elements of the form $d_5(p:0:q)$.  When $q\ne 0$ and $p\ne q$, we also observe that $d_5(p: 0: q)\sim d_5(p:p:q)$, so these two $\P^1$s essentially parameterize the same algebras, with two exceptions.  The algebra $d_5(1: 0:0)$ is isomorphic to $d_5(0:1:0)$ which we saw was isomorphic to $d_{15}(1: 0)$. The algebra $d_5(1: 0:1)$ is not isomorphic to an element outside the family.

The third $\P^1$ consists of elements of the form $d_5(p:p:q)$, and as already noted above, gives elements isomorphic to those of the form
$d_5(p: 0:q)$ except for the elements $d_5(1:1:1)$, which is isomorphic to $d_{14}(0: 1)$, and the element $d_5(1:1: 0)$.

The fourth $\P^1$ consists of elements of the form $d_5(p:q:0)$. We have $d_5(p:q:0)\sim d_5(q:p:0)$, so there is some symmetry on this $\P^1$.
The only nongeneric point is $d_5(1:0:0)=d_5(0:1:0)$, which we have already seen is isomorphic to $d_{15}(1:0)$.

The fifth $\P^1$ consists of elements of the form $d_5(p:q:p)$. When $p\ne \pm q$, we have that $d_5(p:q:p)\sim d_5(q^2:\pm pq:p^2)$, so the
two strata overlap except at the points $(1:1)$ and $(1:-1)$ of the fifth $\P^1$.  Next, we also have a symmetry on this fifth $\P^1$ minus
the divisor $\{(1:1),(1:-1)\}$, given by $d_5(p:q:p)\sim d_5(p:-q:p)$. Note that this is precisely the same divisor we have to exclude when
we are considering the map to the sixth $\P^1$ given by $pr=q^2$, and so there is a corresponding symmetry on that stratum as well.
In fact, the only
elements in that stratum not covered by this map are $d_5(1:1:1)$ and $d_5(1:-1:1)$. The first of these points is isomorphic to $d_{14}(0:1)$,
so the only remaining point is $d_{5}(1: -1: 1)$.

The seventh $\P^1$ consists of elements of the form $d_5(p:q:q)$. This $\P^1$ also collapses to a discrete set of points. If $p\ne0$ and $q\ne0$,
and $p\ne q$, then $d_5(p:q:q)\sim d_5(2:1:1)$.  The cases where $p=0$ or $q=0$ are isomorphic to $d_5(1:0:0)$, which has already been discussed above, and when $p=q\ne0$, we obtain $d_5(1:1:1)$, which has also been discussed above. Thus we only obtain one new case $d_5(2:1:1)$.

In our analysis, we have not taken into account that in addition to the $\lambda$ term, there may be a nontrivial $\psi$ term.  One can check that
for our $\lambda$, the cocycle condition $[\lambda,\psi]=0$ is automatically satisfied. We have $\psi=\psa{45}1c_1+\psa{45}2c_2+\psa{45}3c_3$.
There is a special cases of an algebra where the $\psi$ term doesn't vanish, and we will discuss it later.

\subsubsection{The case $u=0$.} In this case,  the matrix of $B$ reduces to the form $B=\left[\begin{matrix}-rp&r&0\\0&0&r\\0&0&r(q-p)\end{matrix}\right]$.  Now, if $r=0$, we obtain an
 algebra which has a nilradical of dimension higher than 3, so we can reduce to the case when $r=1$. If we take the algebra
$$d_6(p:q)=\psa{24}1+\psa{24}2p+\psa{34}2+\psa{34}3q-\psa{15}1p+\psa{25}1+\psa{35}2+\psa{35}3(q-p),$$
then we obtain a projective family.  However, there is a group of symmetries for $d(p:q)$, which is generated by
\begin{align*}
\sigma(p:q)&=(p-q:p)\\
\tau(p:q)&=(p: p-q).
\end{align*}
These transformations yield a group isomorphic to $\Sigma_3$, with $\sigma$ being of order 3, and $\tau$ of order 2.
However, we also have to remove a divisor from this $\P^1$, consisting of the points $(1:0)$, $(1:1)$ and $(0:1)$, because
the group of symmetries above are not all well-defined on these three points.  Moreover, $d_6(1:1)\sim d_6(1:0)$, so we obtain only two
special algebras $d_6(1:0)$ and $d_6(0:1)$.

We still need to take into account that we may add a $\psi$ term to the $d(p:q)$.  It is easily checked that there are no conditions on the adding of the $\psi$ term, and we obtain one additional algebra:
\begin{align*}
d_7=\psa{24}1+\psa{24}2+\psa{34}2+2\psa{34}3+\psa{35}1+2\psa{35}2+2\psa{35}3+\psa{45}3.
\end{align*}

\subsubsection{Algebras not matching the first matrix form in \refeq{threematrices}}

Now, suppose that no linear combination of the matrices $A$ and $B$ is of the first form in \refeq{threematrices}.  Then either both $A$ and $B$ can
be replaced with linear combinations which are diagonalizable, or we must have a Jordan form of the second type for one of the matrices. Let us
assume it is matrix $A$. Now two commuting matrices can be made simultaneously upper triangular, and we may then conjugate one of them by another
upper triangular matrix to put it in Jordan normal form.

Let the diagonal entries in the upper triangular matrix $B$ be $x$, $y$ and $z$. Now if
$p\ne q$,  then we must have $x=y$, otherwise, a linear combination of $A$ and $B$ would have distinct eigenvalues, and so could be represented by
the first type of matrix in \refeq{threematrices}. Moreover, unless $y=z$ we can assume that $p\ne q$. But in this case, we can add a multiple of
the first matrix to the second to make $y\ne z$. Since the matrix $B$ is assumed to not have a Jordan normal form of the first type, its Jordan
form must be $\left[\begin{smallmatrix}x&0&0\\0&x&1\\0&0&y\end{smallmatrix}\right]$, so that adding a multiple of $A$ to $B$ produces a diagonalizeable
matrix, and since $p\ne q$, we know that $A$ is diagonalizeable.  But we have assumed that we cannot do this, so the alternative is that $p=q$ and
$x=y$.  But now, adding a multiple of one of the matrices to the other gives a strictly upper triangular matrix, which would give an algebra
with nilradical of dimension at least 4.

This means we can assume that both $A$ and $B$ are diagonalizable, and actually are diagonal. Thus we have two matrices of the form
$\diag(p,p,q)$ and $\diag(x,x,y)$.  By adding a multiple of one to the other we can assume that $x=0$, and we also may assume $y\ne0$ because
otherwise the nilradical has too high of a dimension. In fact, we can scale $z=1$, and then by adding a multiple of the first matrix to the
second, we can assume $q=0$ and even $p=1$.  This gives the algebra $d_{8}=\psa{14}1+\psa{24}2+\psa{35}3$.  Finally we ask what happens if we
add a $\psi$ term $\psi=\psa{45}1c_1+\psa{45}2c_2+\psa{45}3c_3$?  It can be shown by direct computation that such an algebra is in fact just
isomorphic to $d_8$.
\subsection{Algebras with 3-dimensional nilradical $\mu=\psa{23}1$}
In this case, the compatibility condition $[\mu,\lambda]=0$ puts some constraints on the matrices $A$ and $B$ of $\lambda$, and taking into account
that we can add a $\mu$-coboundary term to simplify the form of $\lambda$, we first obtain that the $A$ matrix is of the form
$
A=\left[\begin{matrix}a_{22}+a_{33}&0&0\\0&a_{22}&a_{23}\\0&a_{32}&a_{33}\end{matrix}\right]
$, and the $B$ matrix is of the same form with coefficients $b_{ij}$ instead of $a_{ij}$.  However, one can also presume that $a_{23}=0$.
Next, the Maurer-Cartan condition puts additional constraints on the matrices, so we eventually can represent $\lambda$ in the form
$\lambda=2\psa{14}1+\psa{24}2+\psa{34}3+\psa{15}1(p+q)+\psa{25}2p+\psa{25}3+\psa{35}3q$.

Next, we have to take into account the possible $\psi$ term $\psi=\psa{45}1c_1+\psa{45}2c_2+\psa{45}3c_3$. The cocycle condition
$[\lambda,\psi]=0$ forces $c_2=c_3=0$, so we only have to add a single term.

It appears we have generated a whole family of algebras, but it turns out that when $p\ne q$, all of them are isomorphic to
$d_4=\psa{23}1+\psa{14}1+\psa{24}2+\psa{15}1+\psa{35}3$.  When $p=q$ the algebra is isomorphic to $d_9(p:0)$. Thus we only obtain one
new algebra.

This completes the analysis of the elements of the moduli space. Table \ref{Table 3} gives a description of the elements in the moduli
space, along with the dimensions of the cohomology of the elements.  For families, only the generic cohomology is given, as subfamilies
may have cohomology of larger dimension, particularly the cohomology dimension $h^2$ of the cohomology $H^2$ may be larger, meaning that
special subfamilies and elements may have more deformations  than the generic elements.

\begin{table}[htp]
\begin{center}
\begin{tabular}{lrrrrrr}
Type&$H^0$&$H^1$&$H^2$&$H^3$&$H^4$&$H^5$\\
\hline
$d_1=2\psa{34}4-2\psa{35}5+\psa{45}3+\psa{12}1$&0&0&0&0&0&0\\
$d_2=2\psa{34}4-2\psa{35}5+\psa{45}3+\psa{12}1$\\
\hspace{.37in}$-\psa{13}1+\psa{23}2+\psa{24}1+\psa{15}2$&0&1&0&0&1&0\\
$d_3=2\psa{34}4-2\psa{35}5+\psa{45}3$&2&4&2&2&4&2\\
$d_4=\psa{23}1+\psa{14}1+\psa{24}2+\psa{15}1+\psa{35}3$&0&0&0&0&0&0\\
$d_5(p:q:r)=\psa{24}1+\psa{24}2p+\psa{34}2+\psa{34}3q$\\
\hspace{.99in}$+\psa{15}1(q-r)p+\psa{25}1(r-q)+\psa{35}1$\\
\hspace{.99in}$+\psa{35}2r-\psa{35}3r(p-q)$&0&1&2&1&0&0\\
$d_6(p:q)=\psa{24}1+\psa{24}2p+\psa{34}2+\psa{34}3q$\\
\hspace{.71in}$-\psa{15}1p+\psa{25}1+\psa{35}2+\psa{35}3(q-p)$&0&1&2&1&0&0\\
$d_7=\psa{24}1+\psa{24}2+\psa{34}2+2\psa{34}3$\\
\hspace{.34in}$+\psa{35}1+2\psa{35}2+2\psa{35}3+\psa{45}3$&1&2&2&1&0&0\\
$d_8=\psa{14}1+\psa{24}2+\psa{35}3$&0&3&6&3&0&0\\
$d_9(p:q)=\psa{23}1+\psa{34}2+\psa{15}1(2p+q)+\psa{25}2(p+q)$\\
\hspace{.79in}$+\psa{35}3p+\psa{35}4+\psa{45}1+\psa{45}4q$&0&1&1&0&0&0\\
$d_{10}=\psa{23}1+\psa{34}2+3\psa{15}1+2\psa{25}2+\psa{35}3+\psa{45}1+\psa{45}4$&0&2&2&0&0&0\\
$d_{11}=\psa{23}1+\psa{34}2+\psa{15}1+\psa{25}2+\psa{35}4+\psa{45}4$&0&2&4&2&0&0\\
$d_{12}(p:q:r)=\psa{34}2+\psa{15}1p+\psa{15}2+\psa{25}2(q+r)$\\
\hspace{.34in}$+\psa{35}1+\psa{35}3q+\psa{45}3+\psa{45}4r$&0&2&2&0&0&0\\
$d_{13}(p:q)=\psa{34}2+\psa{15}1(p+q)+\psa{25}2(p+q)+\psa{35}1$\\
\hspace{.71in}$+\psa{35}3p+\psa{45}3+\psa{45}4q$&0&3&3&0&0&0\\
$d_{14}(p:q)=\psa{34}2+\psa{15}1p+\psa{15}2+\psa{25}2(p+q)$\\
\hspace{.71in}$+\psa{35}1+\psa{35}3q+\psa{35}4+\psa{45}4p$&0&3&3&0&0&0\\
$d_{15}(p:q)=\psa{34}2+\psa{15}1p+\psa{15}2$\\
\hspace{.71in}$+\psa{25}2(2q)+\psa{35}1+\psa{35}3q+\psa{45}4q$&0&4&4&0&0&0\\
$d_{16}=\psa{34}2+2\psa{15}1+2\psa{25}2+\psa{35}1+\psa{35}3+\psa{45}4$&0&5&5&0&0&0\\
$d_{17}=\psa{34}2+\psa{15}1+\psa{15}2+2\psa{25}2+\psa{35}3+\psa{45}4$&0&6&6&0&0&0\\
$d_{18}=\psa{34}2+\psa{15}1+\psa{25}2+\psa{35}1+\psa{35}4+\psa{45}4$&0&4&8&4&0&0\\
$d_{19}=\psa{34}2+\psa{15}2$&1&11&20&21&15&4\\
$d_{20}(p:q:r:s)=\psa{15}1p+\psa{25}1+\psa{25}2q$\\
\hspace{1.09in}$+\psa{35}2+\psa{35}3r+\psa{45}3+\psa{45}4s$&0&3&3&0&0&0\\
$d_{21}(p:q:r)=\psa{15}1p+\psa{25}2p+\psa{35}2$\\
\hspace{.91in}$+\psa{35}3q+\psa{45}3+\psa{45}4r$&0&5&5&0&0&0\\
$d_{22}(p:q)=\psa{15}1p+\psa{25}1+\psa{25}2q$\\
\hspace{.71in}$+\psa{35}3p+\psa{45}3+\psa{45}4q$&0&7&7&0&0&0\\
$d_{23}(p:q)=\psa{15}1p+\psa{25}2p+\psa{35}3p+\psa{45}3+\psa{45}4q$&0&9&9&0&0&0\\
$d_{24}=\psa{15}1+\psa{25}2+\psa{35}3+\psa{45}4$&0&15&15&0&0&0\\
\hline
\end{tabular}
\end{center}
\caption{\protect Cohomology of Five Dimensional Complex Lie
Algebras}\label{Table 3}
\end{table}

\newpage
\section{Comparison to earlier results}
The first classification of 5-dimensional Real Lie algebras appeared in \cite{mub1}.
We found it difficult to piece out the classification in the paper,
because no list of algebras was included. On the other hand, in \cite{PSWZ} a list of the nondecomposable algebras was given, and the same
list appeared in the recent paper \cite{FGH}. Because the classification is given for Real, rather than Complex, Lie algebras, there is some
duplication if one regards the algebras as complex Lie algebras. For example, the algebras numbered 36 and 37 on the list are two real forms
for the same complex algebra.  We found the list useful in checking our work, and verified that every algebra in the list is covered by our
classification, and no additional algebras were found.  There is mistake in the list in \cite{PSWZ}, because the family numbered $32(u)$, depending
on the parameter $u$, gives algebras which are all isomorphic as real algebras for $u\ne0$, so the inclusion of the parameter $u$ is a mistake.

Our goal in this work is to analyze the structure of the moduli space, in terms of a natural stratification by orbifolds, dictated by deformation
theory.  We have been able to determine such a stratification, where each stratum consists of either a singleton point, a projective orbifold
given by the action of a symmetric group on a $\P^n$ for some $n$ between $1$ and $4$, or a quasi projective orbifold given by an action of
a symmetric group on a quasi projective space obtained by removing a divisor from either $\P^2$ or $\P^1$.  This last phenomenon did not occur
in lower dimension, but the projective orbifolds were present.  The important feature of the stratification is that given an algebra, it will
have smooth deformations along its stratum (assuming it is not a singleton stratum), jump deformations to elements in different strata, and
smooth deformations along that stratum in a nbd of the point to which it jumps.   This gives a natural method of characterizing the moduli space.

The classification in \cite{PSWZ} gives a total of 40 elements, some depending on parameters, and leaves out the decomposable algebras. Our method
of constructing the moduli space gives a total of only 24 strata, consisting of 11 families, of which one is parameterized by an orbifold constructed
from $\P^3$, three from orbifolds constructed from $\P^2$, and seven from orbifolds constructed from $\P^1$, as well as thirteen singletons (not
including the trivial algebra).    In \cite{PSWZ}, there are 18 families
and 22 singletons.  Some of the variation is explained by the fact that we classify complex Lie algebras, and the others classify real Lie algebras,
but even taking this into account, the previous classification does not capture the essence of how the space is glued together. Thus our classification scheme is simpler than the previous ones.

\section{Deformations of the algebras}
In order to determine the deformations of the algebras, we compute the \emph{versal deformation} of an algebra, which is a multi-parameter
deformation that contains the information about all deformations.  The universal infinitesimal deformation of an algebra $d$ is of the form
$$
d^1=d+\delta^it_i,
$$
where the $\delta^i$ give a prebasis of $H^2$, that is, they are 2-cocycles which project to a basis of $H^2$.
The deformation is infinitesimal because $[d^1,d^2]=0\pmod(I_2)$, where $I_2$ is the ideal in $\C[t_1,\mcom]$ generated by all quadratic terms
$t_it_j$, in other words, it is zero up to first order.  The universality means that every \emph{infinitesimal deformation} of $d$ can be
obtained from $d^1$.

The versal deformation is of the form
$$d^\infty=d+\delta^it_i +\gamma^ix_i,$$
where $x_i$ are formal power series in the variables $t_i$ of order at least $2$.
Let $\langle\alpha^i\rangle$ be a prebasis of $H^3$, $\langle\beta^i\rangle$ be a basis of the 3-coboundaries, and $\tau^i$ be a prebasis of the
4-coboundaries, so that jointly, they give a basis of the 3-cochains.   Then for a 2-cochain $d^\infty$, we can express
$$
[d^\infty,d^\infty]=\alpha^ir_i +\beta^is_i+\tau^iu_i.
$$
Then $d^\infty$ is said to be a (mini)versal deformation if $s_i=0$ for all $i$. If we choose $\gamma^i$ so that $D(\gamma^i)=\frac12\beta^i$,
then
$s_i=x_i+ a^{jk}t_jt_k+b^{jk}t_jx_k +c^{jk}x_jx_k$ for some constants $a^{jk}$, $b^{jk}$ and $c^{jk}$. But this means that we can use the
implicit function theorem to show that we can solve the system of equations $s_i=0$ for the $x_i$ as functions of the $t_i$ in some nbd of
the origin in $t,x$ space.  Substituting in the solution in $d^\infty$ gives the versal deformation.  However, we also obtain some relations
$r_i$, which are functions of the $t_i$, with power series having order at least 2.  Finally, it can be shown that the $u_i$ lie in the ideal
generated by the $r_i$. Moreover, the $x_i$ are actually analytic functions of the $t_i$, since the $s_i$ are polynomial in the $t$ and $x$ variables, so $d^\infty$ is a convergent power series. If we take a solution the system of equations $r_i=0$, we obtain an element $d^\infty$ which is an actual
algebra.

We have developed software using the Maplesoft computer algebra system for doing these computations.  Of course, there can be several problems
which occur. First, the system of equations $s_i=0$ may be difficult to solve. Since these are essentially quadratic equations in several variables,
and it is known that the problem of solving quadratic equations is NP complete, meaning that we don't have a polynomial time algorithm for solving
arbitrary systems of quadratic polynomials, it is not surprising that they may be difficult to solve.  Even if we can find a solution, the form
of the solution may be very complicated.  Furthermore, we also have to solve the \emph{relations} $r_i=0$, and this may prove impossible.

For algebras of dimension 4 or less, we were always able to solve the equations $s_i=0$, and surprisingly, we obtained that the $x$ variables
were expressable as rational functions of the $t$ variables.  For 5-dimensional algebras, we have encountered cases where solutions exist that
are not rational, but still are algebraic, and other cases where we could not solve the equations. Thus we cannot give a complete picture
of the deformations, but still we obtain a fairly reasonable approximation of the picture.

We already remarked on the deformations of the algebras $d_1$, $d_2$ and $d_3$. There is a jump deformation from $d_3$ to $d_1$ and the other
two algebras are rigid.  $d_4$ is completely rigid as well since all of its cohomology vanishes.

The algebras $d_5(p:q:r)$ given by the quasi projective orbifold, excluding the seven special $\P^1$s, all have $h^2=2$, so deform only along the family, that is, they have only smooth deformations along the family.
There is one algebra in this group $d_5(3:-3:1)$ whose cohomology dimensions are $0,1,2,3,4,2$ (in other words $h^0=0$, $h^1=1$,\etc, where $h^i$ is the dimension of $H^i$, but since $h^2=2$, it deforms in a generic manner.

Some of the special $\P^1$s have generic deformation behavior.  The subfamily $d_5(p:0:q)$ $d_5(p:q:0)$, $d_5(p:p:q)$ and $d_5(p:q:p)$ as well as
the $\P^1$ given by $d_5(p:q:r)$ such that $pr=q^2$, all have generic cohomology except at some special points. The family $d_5(p:q:q)$ generically
has cohomology numbers $1,3,3,1,0,0$, so has extra deformations, and the two elements in the degenerate subfamily $d_5(0:p:q)$ both reside in other
families, so are not really part of the family $d_5(p:q:r)$.

The family $d_6(p:q)$ has $h^2=3$ generically.  The special point $(1:0)$ has generic cohomology, and so deforms generically, but the special point
$d(0:1)$ has cohomology dimensions $2,7,10,7,2,0$, so has a lot of extra deformations. Also the generic element $d_6(0:0)$ has a lot of extra deformations.

The algebra $d_7$ has $h^2=2$. From its construction, it is evident that it deforms along the family $d_5(p:q:r)$. However, even though
the second cohomology group is of low dimension, we were not able to explicitly construct the versal deformations.  On the other hand,
we computed a (pre)basis $\langle\delta^1,\delta^2\rangle$ of $H^2$, and discovered that the infinitesimal deformations $d^i=d+\delta^it$ are
both actual deformations, and both are smooth deformations along the family $d_5(p:q:r)$ in a nbd of $d_5(1:0:0)$, so that accounts for
all of the deformations of $d_7$.

Now, given the construction of $d_7$, and the fact that $d_5(0:0:0$ is isomorphic to $d_{15}(1:0)$, it
really is evident that $d_7$ and $d_5(1:0:0)$ should be exchanged because they both deform in a nbd of $d_5(0:0:0)$ and the one with the
larger cohomology should not be the one in the family. Since the cohomology numbers of $d_5(1:0:0)$ are $1,5,9,7,2,0$, it is clear which
element would belong in the family. Note that the cohomology $h^2$ of $d_7$ is generic for the family $d_5(p:q:r)$. On the
other hand, the element $d_5(1:0:0)$ doesn't belong to the quasi-projective family, but is a special element, so we did not formally
interchange the elements.

The algebra $d_8$ has $h^2=6$, and it jumps to $d_5(1: 1: 0)$, $d_6(1:1)$ and deforms in  nbds of these points.

The family $d_9(p:q)$ has $h^2=1$ generically, and therefore, generically it deforms only along the family.  There are some special points, whose
cohomology is not generic. The elements   $d_9(1: 0)$, with cohomology numbers $0,1,2,1,0,0$, $d_9(0:1)$, with cohomology numbers $0,1,3,2,0,0$, $d_9(1:-1)$,
with cohomology numbers $0,1,2,1,0,0$, all have extra deformations, which we did not compute.
The three elements  $d_9(3:-2)$, $d_9(2:-3)$ and $d_9(1:-4)$ have cohomology numbers $0,1,1,1,1,0$, and the element
$d_9(3:-4)$, with cohomology $0,1,1,0,1,1$, deform generically, but have different cohomology than the generic case.
The generic element $d_9(0:0)$ has cohomology $0,1,4,7,8,6,2$, and thus a lot of deformations. It is nilpotent, and also isomorphic to
$d_{12}(0:0:0)$.

The algebra $d_{10}$ has $h^2=2$, and there is a jump deformation to $d_9(1:1)$ and it also deforms in a nbd of this point.

The algebra $d_{11}$ has $h^2=$ and it is \emph{hard to compute the deformations.}

The family $d_{12}(p:q:r)$ generically has $h^2=2$, so generically, algebras in this family have deformations only along the family.
There are some special $P^1$s, whose cohomology (and deformations) are not generic.

Generically, elements of the family $d_{13}(p:q)$ have $h^2=3$, and generically deform along the family, jump to $d_{12}(p,q,p+q)$ and deform in a nbd of that algebra.  There are also a lot of special cases, but we won't study them here.

Generically, elements of the family $d_{14}(p:q)$ have $h^2=3$, deform along the family, jump to $d_{12}(p:p:q)$ and deform in a nbd of this algebra.
There are several special values of the parameters that have nongeneric deformation patterns.

Generically, elements of the family $d_{15}(p:q)$ have $h^2=4$, deform along the family, jump to $d{12}(p:q:q)$ and deform in a nbd of this algebra.

The algebra $d_{16}$ has $h^2=5$, it jumps to $d_{12}(2:1:1)$, $d_{13}(1:1)$ and $d_{15}(2: 1)$ as well as deforming in nbds of these points.

The algebra $d_{17}$ has $h^2=8$, it jumps to $d_{12}(1:1:1)$, $d_{14}(1:1)$ and $d_{15}(1:1)$ as well as deforming in nbds of these points.

The algebra $d_{18}$ has $h^2=8$ and we were unable to compute the versal deformation. We obtained some information on the deformations, by
computing deformations given by using several, but not all of the cocycles at the same time. These are not really the versal deformation, but are
multi-parameter deformations which do reveal a lot about the general deformation picture.
There are jump deformations to $d_6(x: y)$, $d_8$, $d_{11}$, $d_{12}(2: 1: 1)$, $d_{13}(1: 1)$, $d_{14}(1: 0)$ and $d_{15}(2: 1)$ and deforms in a nbd of $d_{12}(2: 1: 1)$, $d_{13}(1: 1)$, $d_{14}(1: 0)$, and $d_{15}(2: 1)$.

The algebra $d_{19}$ has $h^2=20$, and so it is surprising that we were able to compute the versal deformation easily.  However $h^3=21$, so there
are 21 relations in 20 variables, which we were not able to solve.  What we do in this case is reduce the number of variables by setting some of them
equal to zero, and then solving. This gives a partial deformation, in the sense that the solutions give some of the deformations, but not necessarily all of them.  We then apply the same method with other variables being set equal to zero, and so we obtain more of the solutions. With a little
luck, we may be able to obtain all of the solutions. In particular, we found jump deformations to $d_9(0:0)$, $d_{12}(0: 0: 0)$, $d_14(0:0)$, and $d_{15}(0:0)$, which means it also must jump to $d_5(0:0:0)$ because it is isomorphic to $d_{14}(0:0)$.  Now, the generic element in a family
has jump deformations to every element in the family, so $d_{19}$ must have jump deformations to every element of the form
$d_5(x:y:z)
$, $d_9(x:y)$, $d_{12}(x:y:z)$, $d_{14}(x:y)$, and $d_{15}(x:y)$. In particular, we know which other nilpotent algebras $d_{19}$ deforms to.

The family $d_{20}(p:q:r:s)$ has $h^2=3$, which means that generically, the deformations lie only along the family. There are some special
subfamilies whose cohomology is not generic and which have more deformations, but there are a lot of them, and we will not include this information.

The family $d_{21}(p:q:r)$ has $h^2=5$, and generically, there is a jump deformation to $d_{20}(p:p:q:r)$, and deformations in a nbd of this point,
as well as deformations in a nbd of $d_{21}(p:q:r)$. Again there are some special subfamilies, which we do not describe here.
The $4\times4$ matrix $A$ which gives the versal deformation for a generic element is given by
$$
A=\left[ \begin {array}{cccc} p&t_{{4}}&0&0\\\noalign{\medskip}0&p+t_{{2}}&1&0\\\noalign{\medskip}t_{{5}}&t_{{3}}&q&1\\\noalign{\medskip}0&t_
{{1}}&0&r\end {array} \right]
.
$$

The family $d_{22}(p:q)$ has $h^2=7$, and generically, there are jump deformations to $d_{21}(p:q:q)$, $d_{21}(q:p:p)$ and $d_{20}(p:p:q:q)$, as well
as deformations in nbds of these points.
The $4\times4$ matrix $A$ which gives the versal deformation for a generic element is given by
$$
A=\left[ \begin {array}{cccc} p&1&t_{{7}}&0\\\noalign{\medskip}0&q+t_{{5}}&0&t_{{3}}\\\noalign{\medskip}0&0&p+t_{{1}}&1\\\noalign{\medskip}t_
{{4}}&t_{{6}}&t_{{2}}&q\end {array} \right]
.
$$

The family $d_{23}(p:q)$ has $h^2=9$, and generically, there are jump deformations to $d_{21}(p:p:q)$ and $d_{20}(p:p:p:q)$, and it also deforms in nbds of these points.
The $4\times4$ matrix $A$ which gives the versal deformation for a generic element is given by
$$
A= \left[ \begin {array}{cccc} p+t_{{2}}&t_{{4}}&t_{{8}}&0
\\\noalign{\medskip}t_{{1}}&p+t_{{6}}&t_{{9}}&0\\\noalign{\medskip}0&0
&p&1\\\noalign{\medskip}t_{{3}}&t_{{5}}&t_{{7}}&q\end {array} \right]
.
$$

The algebra $d_{24}$, with $h^2=15$, has jump deformations to $d_{23}(1:1)$, $d_{22}(1:1)$, $d_{21}(1:1:1)$ and $d_{20}(1:1:1:1)$ as well as deformations in nbds of these points.  Note that the $4\times 4$ matrix which determines an extension of the trivial 1-dimensional algebra
by the trivial 4-dimensional algebra has 16 entries, and by the projectivity of the space of such algebras, we could have a maximum of 15
independent cocycles, which is exactly realized by this algebra. Interestingly, the algebra only deforms within this subclass of algebras.
The versal deformation can be given by the $4\times 4$ matrix $A$ below.

$$A=\left[ \begin {array}{cccc} 1+t_{{14}}&t_{{7}}&t_{{5}}&t_{{1}}\\\noalign{\medskip}t_{{11}}&1+t_{{8}}&t_{{6}}&t_{{15}}
\\\noalign{\medskip}t_{{12}}&t_{{9}}&1+t_{{2}}&t_{{4}}
\\\noalign{\medskip}t_{{13}}&t_{{10}}&t_{{3}}&1\end {array} \right].
$$
It should be evident from the form of this matrix that $d_{24}$ deforms to all of the families $d_{20}-d_{23}$.
\section{Nilpotent Algebras}
There are 11 families, $d_5(p:q:r)$, $d_6(p:q)$, $d_{12}(p:q:r)$, $d_{13}(p:q)$, $d_{14}(p:q)$, $d_{15}(p:q)$, $d_{20}(p:q:r:s)$, $d_{21}(p:q:r)$,
$d_{22}(p:q)$ and $d_{24}(p:q)$. The generic element in each of these families is nilpotent.  In addition, the algebra $d_{19}$, which is not a member of a family, is also nilpotent. This would give 12 nilpotent algebras, except that some of these algebras are isomorphic. The isomorphisms are
\begin{align*}
d_5(0:0:0)&\sim d_{14}(0: 0)\sim d_{15}(0:0)\\
d_6(0: 0)&\sim d_{21}(0:0:0)\\
d_9(0:0)&\sim d_{12}(0:0:0)
\end{align*}
Thus we eliminate 4 of the algebras, as duplicates, which leaves 8 algebras.  The trivial algebra is not in this list of algebras, and it is
also abelian, so in the end, we have a total of 9 nonisomorphic nilpotent complex Lie algebras.  The algebras $d_{21}(0:0:0)=\psa{35}2+\psa{45}3$,
$d_{23}=\psa{45}3$ and the trivial algebra, so there are 6 indecomposable algebras.

It would be very difficult to calculate the versal deformations of these algebras. However, we give a table of their cohomology, which should give
some idea of how complicated their deformations are.
\begin{table}[htp]
\begin{center}
\begin{tabular}{lrrrrrr}
Type&$H^0$&$H^1$&$H^2$&$H^3$&$H^4$&$H^5$\\
\hline
$\psa{24}1+\psa{34}2+\psa{35}1$&1&6&13&15&10&3\\
$\psa{35}2+{45}3$&2&8&14&15&10&3\\
$\psa{34}2+\psa{15}2+\psa{35}1+\psa{45}3$&1&4&7&8&6&2\\
$\psa{34}2+\psa{35}1+\psa{45}3$&2&7&9&9&7&2\\
$\psa{34}2+\psa{15}2$&1&11&20&21&15&4\\
$\psa{25}1+\psa{35}2+\psa{45}3$&1&5&8&8&6&2\\
$\psa{25}1+\psa{45}3$&2&10&19&20&12&3\\
$\psa{45}3$&3&14&28&13&17&4\\
\hline
\end{tabular}
\end{center}
\caption{\protect Cohomology of Nilpotent Five Dimensional Complex Lie
Algebras}\label{Table 4}
\end{table}

We don't give the cohomology for the trivial algebra, because every cochain is a cocycle, so we have $h^n=5^{n+1}$, and of course, the trivial algebra
deforms to every other algebra.

The important observation to make here is  that all but one of our nilpotent algebras arise in a family of solvable algebras, a fact that wasn't
observable in the previous classifications.
\section{Conclusions and Questions}

In the moduli spaces of Lie algebras of dimension 3, 4 and 5, we found a simple description of a stratification of the moduli space which corresponds
in a clear manner to the deformations of the algebras. We were able to give a numerical ordering of the algebras such that an algebra only deforms to
algebras whose number is less than or equal to the number of the algebra, with one exception, that the generic element in a family might not follow this restriction.  In some sense, this means  the generic element in a family does not precisely belong to the family.

We also saw that every family of solvable Lie algebras contained a nilpotent algebra, and moreover, it always was the generic element in the family.
On the other hand, there is a nilpotent 5-dimensional Lie algebra which does not belong to a family. We conjecture that families of nilpotent algebras
are special subfamilies of families of solvable Lie algebras, and that every family of solvable Lie algebras will contain at least one nilpotent
algebra.

Based on the results we obtained about moduli spaces of low dimensional complex Lie and associative algebras, the authors came to a conjecture
that these moduli spaces have a stratification by projective orbifolds, ``glued'' together by jump deformations. This conjecture has been shown
to be true for Lie and associative algebras of dimension up to 4, and for non nilpotent associative algebras of dimension 5. Here, we showed that
the conjecture is still true for complex Lie algebras of dimension 5, but we had to allow quasi projective orbifolds, so the picture is a bit
more complicated. It is possible that the quasi projective orbifolds could be replaced with projective orbifolds, but we don't know if this is
possible.  We are reaching the limit of the computational power of our software, and a proof of our conjecture would need a more sophisticated
understanding of the moduli spaces.

\bibliographystyle{amsplain}
\providecommand{\bysame}{\leavevmode\hbox to3em{\hrulefill}\thinspace}
\providecommand{\MR}{\relax\ifhmode\unskip\space\fi MR }
\providecommand{\MRhref}[2]{%
  \href{http://www.ams.org/mathscinet-getitem?mr=#1}{#2}
}
\providecommand{\href}[2]{#2}

\end{document}